\documentclass[journal,technote]{IEEEtran}
\pdfoutput=1
\usepackage{graphicx}
\usepackage{subcaption}
\usepackage{amsmath}
\usepackage{epstopdf}
\usepackage{cite}
\usepackage{amssymb}
\usepackage{bm}
\usepackage{color}
\usepackage{algorithm}
\usepackage{algorithmicx}
\usepackage{algpseudocode}
\usepackage{multirow}
\usepackage{makecell}
 \usepackage{url}

\newtheorem{theorem}{\bf{\textit{Theorem}}}
\newtheorem{lemma}{\bf{\textit{Lemma}}}
\newtheorem{remark}{\bf\textit{Remark}}
\newtheorem{assumption}{\bf\textit{Assumption}}

\renewcommand{\theequation}{\arabic{equation}}
\renewcommand{\thefigure}{\arabic{figure}}
\newtheorem{problem}{\bf\textit{Problem}}

\renewcommand{\thetable}{\arabic{table}}

\newcommand*{\QEDB}{\hfill\ensuremath{\square}}
\newcommand*{\QEDD}{\hfill\ensuremath{\bullet}}

\algnewcommand{\Initialize}[1]{%
  \State \textbf{Initialize:}
  \Statex \hspace*{\algorithmicindent}\parbox[t]{1.0\linewidth}{\raggedright #1}
}

\renewcommand\baselinestretch{1}
\ifCLASSINFOpdf
\else
\fi
\begin{document}
%
\title{Robust Output Regulation and Reinforcement Learning-based Output Tracking Design for Unknown Linear Discrete-Time Systems} 

%
%

\author{Ci Chen, Lihua Xie, Yi Jiang, Kan Xie, and Shengli Xie

\thanks{This work was supported in part by the Wallenberg-NTU Presidential Postdoctoral Fellowship and in part by National Natural Science Foundation of China under Grants 61703112, 61973087, and 61727810. (Corresponding author: Lihua Xie).}
\thanks{C. Chen and L. Xie are with School of Electrical and Electronic Engineering, Nanyang Technological University, Singapore 639798. (e-mail: ci.chen@control.lth.se, ELHXIE@ntu.edu.sg). }
\thanks{Y. Jiang is with State Key Laboratory of Synthetical Automation for Process Industries, Northeastern University, Shenyang, China. (e-mail: yijiang.ezhou@gmail.com).}
\thanks{K. Xie and S. Xie are with School of Automation, Guangdong University of Technology, Guangdong Key Laboratory of IoT Information Technology, Guangzhou, 510006 China (e-mail: shlxie@gdut.edu.cn).}
}

\maketitle

\begin{abstract}
In this paper, we investigate the optimal output tracking problem for linear discrete-time systems with unknown dynamics using reinforcement learning and robust output regulation theory. This output tracking problem only allows to utilize the outputs of the reference system and the controlled system, rather than their states, and differs from most existing tracking results that depend on the state of the system.
The optimal tracking problem is formulated into a linear quadratic regulation problem by proposing a family of dynamic discrete-time controllers. Then, it is shown that solving the output tracking problem is equivalent to solving output regulation equations, whose solution, however, requires the knowledge of the complete and accurate system dynamics. To remove such a requirement, an off-policy reinforcement learning algorithm is proposed
using only the measured output data along the trajectories of the system and the reference output.
By introducing re-expression error and analyzing the rank condition of the parameterization matrix, we ensure the uniqueness of the proposed RL based optimal control via output feedback.
\end{abstract}


%
\IEEEpeerreviewmaketitle

\section{Introduction}
\IEEEPARstart{O}{utput tracking}, whose objective is to make the system output 
follow a desired reference trajectory,
is a fundamental research topic of practical importance (see examples in \cite{stevens2015aircraft}). One systematic way to approach the output tracking is to transform it into an output regulation problem, whose solution and  corresponding  control design date back to \cite{francis1977linear}, where both the closed-loop stability and the asymptotic tracking of even an unbounded reference are achieved. Though elegant, such a solution is built on the knowledge of the complete and accurate system dynamics, making the output tracking design model dependent.

Reinforcement learning (RL) features making sequential decisions  through interactions between the agent's actions and unknown environment \cite{sutton1998introduction,lewis2012optimal}. RL algorithms have been applied in the control field to solve optimal control problems for both
{discrete-time (DT) (see, e.g., \cite{lewis2012optimal,zhang2012adaptive})} and {continuous-time (CT) systems (see, e.g.,  \cite{lewis2012optimal,zhang2012adaptive,vrabie2013optimal,Reza2016CYbOutput,jiang2017robust,gao2016tacadaptive,kamalapurkar2018RLBook,Ci2019TACRL,Ci2019OPFBRL1,Ci2020AutoRL,Jiang2020book})} without any knowledge of the system dynamics. RL-based methods in {Chapter 11 of \cite{lewis2012optimal}} were used in \cite{kiumarsi2015actor,kiumarsi2014Auto} to handle the optimal tracking control of linear and nonlinear DT systems. However, most RL algorithms are classified as on-policy learning, which assumes that the behaviour policy for generating the data for learning is the same as the target policy. Off-policy RL differs from the on-policy learning in its separating the target and behaviour policies. Off-policy RL for linear
CT systems with unknown dynamics was given
in {Chapter 2 of \cite{jiang2017robust}} to solve the optimal regulation problem. A solution to the zero-sum game problem for the regulation of DT systems was given in  \cite{kiumarsi2017h} using off-policy RL. Together with RL and the output regulation theory, \cite{jiang2019optimal} and \cite{gao2018DTPcopy} gave tracking controllers for single and multiple DT systems and established the asymptotic stability of the tracking error, which is the DT version of \cite{gao2016tacadaptive}. Note that one essential assumption for \cite{jiang2019optimal,gao2018DTPcopy}, as well as \cite{gao2016tacadaptive}, is that not only the outputs of the reference and
the controlled system but also their states are required in {the optimal control design process}. Therefore, a key challenge to be addressed is how to achieve the optimal output tracking of DT systems without using the reference state.

In contrast to full states feedback, output feedback utilizes the system output and adds extra flexibilities to control design. As for the output-feedback control, \cite{LewisCYB2011}  utilized the state reconstruction to form an output-feedback RL for regulating the states of DT systems. Based on \cite{LewisCYB2011}, off-policy $H_\infty$ control of DT systems was given based on input and output data in the literature such as \cite{Jiangyi2018Hinfty,rizvi2018Auto}.  The optimal tracking control using output-feedback RL  algorithm  was given in \cite{Bahare2015CYB}. \cite{Rizvi2018TNNLS} gave a parameterization of the state of a DT system in terms of the input and output data for the output feedback LQR by $Q$-learning. Based on \cite{Rizvi2018TNNLS}, \cite{rizvi2018output} considered an output-feedback tracking control for DT systems with the assumption that the reference state is available during the learning and control processes. {Utilizing the state reconstruction in \cite{LewisCYB2011}, \cite{gao2016OPFBDC,gao2018adaptive} considered tracking-based RL algorithms for disturbance rejection with the augmented system being observable. It turns out that, in the state parametrization, the full row rank of the parametrization matrix is required, see \cite{gao2014adaptive,Reza2016CYbOutput,Ci2019OPFBRL1,rizvi2020outputAuto,rizvi2020tac,rizvi2020book} for example.  {In \cite{rizvi2020outputAuto,rizvi2020tac,rizvi2020book}, sufficient conditions were established that guarantee such a rank condition in both the CT and DT settings.} Within the framework of output-feedback RL, most of the existing output tracking results for DT systems are based on the linear output regulation theory, which,  however, remains challenging for its extension to robust output regulation.}

Motivated by the analysis above, we investigate the output tracking problem of DT systems using RL and {robust output regulation theory} based on the reference output. We use real-time data collected along the trajectories of the DT system and propose an RL algorithm that ensures not only the stability of the closed-loop system, but also the predefined system performance. The contributions of this paper are given as follows.
\begin{enumerate}
    \item {Compared to the linear output regulation based RL controllers for tracking DT systems}  \cite{jiang2019optimal,gao2018DTPcopy,gao2016OPFBDC,rizvi2018output,gao2018adaptive}, we amend the {standard robust output regulation theory} to propose a new family of {robust} DT controllers, based on which an optimal output tracking problem is formulated.
    \item   We invoke the CT work in \cite{Ci2019OPFBRL1} to propose a criterion for guaranteeing the parameterization matrix in the DT state reconstruction to be full row rank, which leads to the uniqueness of the data-based optimal output-feedback DT controller. Inspired by \cite{Ci2019OPFBRL1}, we establish
        the rank of the parameterization matrix for DT systems, and based on which we show that the controllability is a sufficient condition for ensuring the full row rank of the parameterization matrix for DT systems.
     \item We do not need to compute analytic solutions to the output regulation equations. We take the DT state re-expression error into account during the optimal control learning. We add a phase of model-free pre-collection into the off-policy RL for reducing the adverse
impact of the re-expression error.
\end{enumerate}

\textit{Notation:} Throughout this paper, given any square matrix
$\mathcal{M}$, the notation $\lambda (\mathcal{M})$ indicates the spectrum of $\mathcal{M}$, $\rho(\mathcal{M})$ is its spectral radius, $\det(\mathcal{M})$ denotes its determinant, ${\rm{adj}}(\mathcal{M})$ is its adjugate, $\mathcal{M}>0$ $(\mathcal{M}\geq 0)$ means that the matrix is positive definite (positive semi-definite), ${\rm{vec}}(\mathcal{M}) = [\mathcal{M}_1^T ,\mathcal{M}_2^T ,\cdots
,\mathcal{M}_n^T ]^T$ is a vector with $\mathcal{M}_i$ being the $i$th column of $\mathcal{M}$, and ${{\rm{vecs}}(\mathcal{W})} = [\mathcal{W}_{11} ,2\mathcal{W}_{12} ,...,2\mathcal{W}_{1n} ,\mathcal{W}_{22} ,2\mathcal{W}_{23}
,...,2\mathcal{W}_{n - 1,n} ,\mathcal{W}_{n,n} ]^T \in \mathbb{R}^{\frac{1}{2}n(n+1)}$ with $\mathcal{W}_{ij}$ being the entry at the $i$th row and $j$th column of $\mathcal{W}$.
Given any vector $\mathcal{T} \in \mathbb{R}^n$, ${{\rm{vecv}}(\mathcal{T})} = [\mathcal{T}_1^2 ,\mathcal{T}_1 \mathcal{T}_2 ,...,\mathcal{T}_1 \mathcal{T}_n ,\mathcal{T}_2^2,\mathcal{T}_2 \mathcal{T}_3 ,...,\mathcal{T}_{n - 1} \mathcal{T}_n ,\mathcal{T}_n^2 ]^T \in \mathbb{R}^{\frac{1}{2}n(n+1)}$. Let the controllability matrix $\mathcal{C}(A_1, A_2)=[A_2, {A_1} A_2$, ${A_1^2} A_2, \cdots, A_1^{n-1} A_2]$ with
the dimension of $A_1$ being $n\times n$. The notation $ \otimes $
denotes the Kronecker product. The notations $0$ and $I$, respectively, indicate a zero matrix and an identity matrix with appropriate dimensions.

\section{Problem Formulation and Preliminaries}\label{Sect2}

\subsection{Problem Formulation}
We aim to study a class of DT dynamical systems modeled by
\begin{align}
x(k+1) &= Ax(k) +Bu(k),\label{eqsystem7}\\
y(k) &= Cx(k),\label{eqsystem71}
\end{align}
where $x(k)\in\mathbb R^{r_n}$, $u(k)\in\mathbb R^{r_m}$, and $y(k)\in\mathbb R^{r_p}$ denote the system state, input, and output, respectively, and,  $A\in\mathbb R^{r_n\times r_n}$, $B\in\mathbb R^{r_n\times r_m}$, and $C\in\mathbb R^{r_p\times r_n}$, involved in the system dynamics (\ref{eqsystem7}) and (\ref{eqsystem71}), are unknown constant matrices. The control objective is to drive $y(k)$ in (\ref{eqsystem71}) to follow a reference signal $y_d(k)$ given by
\begin{align}
 x_d(k+1) &= {S}x_d(k),\label{eqsystem5}\\
y_d(k) &=Rx_d(k),\label{eqsystem72}
\end{align}
where $x_d(k)\in\mathbb R^{r_{q_m}}$ and $y_d(k)\in\mathbb R^{r_p}$ respectively denote the reference state and output.
Similar to the $(A, B, C)$ in (\ref{eqsystem7}) and (\ref{eqsystem71}), ${S}\in\mathbb R^{r_{q_m}\times r_{q_m}}$ and ${R}\in\mathbb R^{r_p\times r_{q_m}}$ are unknown constant matrices. The output tracking error between (\ref{eqsystem71}) and (\ref{eqsystem72}) is defined as
  \begin{align}
        y_e(k)&=y(k)-y_d(k).\label{eqsystem8}
 \end{align}

For the purpose of control design, some assumptions on the dynamics of DT system  (\ref{eqsystem7})--(\ref{eqsystem72}) are made as follows.
\begin{assumption}\label{assumption3}
For system (\ref{eqsystem7})--(\ref{eqsystem71}), $(A, B)$ is controllable and $(A, C)$ is observable.
\end{assumption}
\begin{assumption}\label{assumption2}
The eigenvalues of the matrix $S$ are on or outside the unit circle.
\end{assumption}
{
\begin{assumption}\label{assumption1}
The minimal polynomial of $S$ is known.
\end{assumption}}
\begin{assumption}\label{assumption4}
{\rm{rank}}$\left(\left[\begin{array}{cc}
 A-\lambda_i I  &\ B \\
    C  &\ 0
     \end{array}\right]\right)=r_n+r_p$, $\forall \lambda_i\in\lambda({S})$.
\end{assumption}

\textit{Assumption \ref{assumption3}} is standard in the optimal control\cite{lewis2012optimal}.
 \textit{Assumption \ref{assumption2}} is made to rule out the trivial  case wherein the matrix $S$ in (\ref{eqsystem5}) is Schur {(see Chapter 1 of \cite{huang2004nonlinear})}, namely,  the reference state is asymptotically stable. {\textit{Assumptions \ref{assumption2}--\ref{assumption4}} are standard in the output regulation literature.}

Based on the system descriptions above, we focus on formulating an optimal tracking control problem that is solvable. To do this, we extend the work of CT systems \cite{Ci2019OPFBRL1} and construct a DT control protocol
\begin{subequations}
 \begin{align}
 z(k+1) &= Fz(k) - Gy_e(k), \hfill \\
u(k)&=-Kx(k)-Hz(k)-Tz(k), \hfill
\end{align}
\label{eqcontroller1}
\end{subequations}
\noindent where {$z(k)\in\mathbb R^{r_pr_{q_m}}$ is a dynamical signal driven by output error $y_e(k)$}, $(F, G)$ is {an} $r_p$-copy internal model of $S$ \footnote{The design of $(F, G)$ is available under \textit{Assumption \ref{assumption1}}; see Chapter 1 of \cite{huang2004nonlinear}.} {with $F\in\mathbb R^{r_pr_{q_m}\times r_pr_{q_m}}$ and $G\in\mathbb R^{r_pr_{q_m}\times r_p}$},
  ${T\in\mathbb R^{r_m\times r_pr_{q_m}}}$ is a newly proposed {feedforward} gain matrix compared to the standard DT controller by the output regulation theory \cite{huang2004nonlinear},
 and {$K\in\mathbb R^{r_m\times r_n}$ and $H\in\mathbb R^{r_m\times r_pr_{q_m}}$} are gain matrices for solving the tracking problem to be specified later.

Substituting the dynamical controller (\ref{eqcontroller1}) into the DT system (\ref{eqsystem7}) yields
\begin{subequations}\label{eq15}
\begin{align}
   {\bar x}(k+1) &= {\bar A}(k){\bar x}(k)   + {{\bar G}}R{x_d}(k),\\
  {y_e}(k) &= {{\bar C}}{{\bar x}}(k) - R{x_d}(k),
\end{align}
\end{subequations}
where {${\bar x(k)}=[{x}(k)^T, {z}(k)^T]^T\in\mathbb R^{n_z}$ with $n_z=r_n+r_pr_{q_m}$}, and the system matrices {$\bar A\in\mathbb R^{n_z\times n_z}$, $\bar G\in\mathbb R^{n_z\times r_p}$, and $\bar C\in\mathbb R^{r_p\times n_z}$} are, respectively, given as
$\bar A=\left[ {\begin{array}{*{20}{c}}
  {{A} - {B}{K}}& {-{B}{H}}-BT \\
  { - {G}{C}}& {{F}}
\end{array}} \right]$,    $\bar G=\left[ {\begin{array}{*{20}{c}}
  0 \\
  {{G}}
\end{array}} \right],\ \ \mbox{and}\ \ \bar C=\left[ {\begin{array}{*{20}{c}}
  {{C}}, & 0
\end{array}} \right]$.
Note that
\begin{align}\label{cceq15001}
{\bar A}=\underline{A}-\bar B[K, H],
\end{align}
 where $\underline{A}=\left[ {\begin{array}{cc}
  {A}&\ -BT\\
   -{G}{C}&\ F
\end{array}} \right]{\in\mathbb R^{n_z\times n_z}}$ and $\bar B=\left[ {\begin{array}{c}
  B\\
  0
\end{array}} \right]{\in\mathbb R^{n_z\times r_m}}$. By \cite{Ci2019OPFBRL1}, the following properties of the system dynamics hold.
\begin{lemma}\label{lemma0}
Under \textit{Assumption \ref{assumption4}} and when $(F, G)$ incorporates {an} $r_p$-copy internal model of $S$, we have that
\begin{enumerate}
  \item given any matrix $T$, if $(A, B)$ is stabilizable (controllable), then $(\underline{A}, \bar B)$ is stabilizable (controllable).
  \item given the matrix $T$ such that $(F, T)$ is observable with $r_p\geq r_m$, if $(A, C)$ is detectable (observable), then the pair $(\underline{A}, \bar C)$ is detectable (observable). In addition, if $(A, B)$ is stabilizable (controllable), then $(\underline{A}, \bar B)$ is stabilizable (controllable).
\end{enumerate}
\end{lemma}

  Given a Schur matrix $\bar A$, under  \textit{Assumptions \ref{assumption3}}--\textit{\ref{assumption2}}, the  equations ${{\bar A}}{X} + {{\bar G}}R = {X}S$ and ${{\bar C}}{X} = R$ have a common solution $X$ \cite{huang2004nonlinear}. Based on the matrix ${X}$, let
\begin{align}\label{eqerror1pre}
{e}(k) = {\bar x}(k) - {X}{x_d}(k),
\end{align}
 which leads to
\begin{align}
  e(k+1) =&  {{\bar A}}{{e(k)}},
  \label{eq18pre}
\end{align}
where (\ref{eq15}) and (\ref{eqerror1pre}) are used. Substituting  (\ref{cceq15001}) into (\ref{eq18pre}) yields
\begin{subequations}\label{eq1901pre}
\begin{align}
  e(k+1) &= \underline{A}{{e(k)}}+\bar B u_e(k),   \label{eq1901pre01} \\
  {y_e}(k) &= {{\bar C}}{{e}}(k), \label{eq1901pre02}\\
   u_e(k) &=-\bar K {e}(k),\label{eq1901pre03}
\end{align}
\end{subequations}
where $\bar K=[K, H]\in\mathbb R^{r_m\times (r_m+r_pr_{q_m})}$ is the feedback gain matrix with $K$ and $H$ being from (\ref{eqcontroller1}). Note that if $e(k)$ decays to zero, so does ${y_e}(k)$. Here, the matrix $\bar K$ in (\ref{eq1901pre03}) is designed based on the following optimization problem.
 \begin{problem} \label{def2}
\begin{align}
&\min_{u_e}\sum_{i=k}^{\infty}({y_e^T}(i)Q{y_e}(i)+u_e^T(i)\bar R u_e(i)) \label{22eqopt1}\\
&\rm{subject\ \ to\ \ } (\ref{eq1901pre}) \nonumber
\end{align}
\rm{with} $Q> 0$ \rm{and} $\bar R>0$.
\end{problem}

\subsection{Preliminaries on Optimal Control}
  The optimal tracking problem is solved if  $\bar K$ in  (\ref{eq1901pre03}) is designed such that $e(k)$ in (\ref{eq1901pre02}) is stabilized to zero and,  meanwhile, the performance index in (\ref{22eqopt1}) is minimized.

   The optimal feedback gain matrix that solves the optimization problem (\ref{22eqopt1}) is labelled as $\bar K^*$ and is given as follows {(see Chapter 2 of  \cite{lewis2012optimal})},
\begin{align}
\bar K^*=(\bar R + \bar B^TP^*\bar B)^{-1}\bar B^TP^*\underline{A}, \label{eq1601}
\end{align}
where $P^*$ satisfies the DT algebraic Riccati equation (ARE)
\begin{align}\label{eq1701}
&\underline{A}^T P^* \bar B (\bar R+\bar B^TP^*\bar B)^{-1}\bar B^TP^*\underline{A}
\nonumber\\&= \bar C^T\bar Q\bar C + \underline{A}^TP^*\underline{A}-P^*
\end{align}
with $\bar Q=\mbox{diag}\{Q, 0\}$
and
 assuming that the stabilizability condition of $(\underline{A}, \bar B)$  and the observability condition of $(\underline{A}, {\bar C})$ hold.  Note that the detectability condition of $(\underline{A}, {\bar C})$ may not hold through the design of  the output regulation-based standard DT controller (see {Chapter 1 of \cite{huang2004nonlinear}}).

 As for solving the DT ARE (\ref{eq1701}),  a model-based algorithm with the knowledge of $\underline{A}$ and $\bar B$ was given in \cite{HewerDT1971}, and is recalled below.
\begin{lemma}\label{lemma1}(\cite{HewerDT1971})
Let $\bar K^0$ be a stabilizing gain matrix such that $\underline{A}-\bar B{\bar K}^0$ is Schur. Solve $P^j$ from
\begin{align}\label{eq18}
\hspace*{-0.1cm}P^j=\bar Q + (\bar K^j)^T\bar R \bar K^j + (\underline{A}-\bar B \bar K^j)^TP^j(\underline{A}-\bar B \bar K^j).
\end{align}
Update the policy as
\begin{align}\label{eq19}
\bar K^{j+1}=(\bar R + \bar B^TP^j\bar B)^{-1}\bar B^TP^j\underline{A}.
\end{align}
Then,
\begin{enumerate}
  \item $P^*\leq P^{j+1}\leq P^j$ for each $j=1, 2, \ldots$;
  \item $\lim_{j\rightarrow\infty}\bar K^j=\bar K^*$, $\lim_{j\rightarrow\infty}P^j=P^*$.
\end{enumerate}
\end{lemma}

\section{Output-Feedback  RL for Optimal Output Tracking Control of DT Systems}\label{Sect4}

This section is to solve \textit{Problem \ref{def2}} using the data collected along the trajectories of the controlled system and the reference output in the absence of any knowledge of the system dynamics. To achieve this, we use the following behavior policy to excite the DT system (\ref{eqsystem7})
{\begin{align}
\bar u(k)=-{\bar K}^0r(k)+\xi(k)-T\bar z(k),
\label{ceqcontroller1}
\end{align}
where ${r}(k)={\left[ {{x}^T(k),{\bar z}^T(k)} \right]^T}\in\mathbb R^{n_z}$, ${\bar K}^0$ is an initial stabilizing gain, $\xi(k)$ is an exploration noise, and $\bar z(k)\in\mathbb R^{r_pr_{q_m}}$ is an alternative dynamical signal driven by the output error $y_e(k)$ as defined by
${\bar{z}}(k+1) = F\bar {z}(k) - Gy(k) +  G\vartheta(k)$
 with
 {$F$ and $G$} being given from (\ref{eqcontroller1}) and $\vartheta(k)$ being either the exploration noise or the reference's output $y_d(k)$. Here, $\bar z(k)$ is defined to differ from $z(k)$ in (\ref{eqcontroller1}). We will specify $\vartheta(k)$ in our design later.}
  Note that, if $\vartheta(k)=y_d(k)$ and $\xi(k)=0$, then (\ref{ceqcontroller1}) is equivalent to (\ref{eqcontroller1}). It follows from (\ref{eqsystem7})--(\ref{eqsystem72}) and (\ref{ceqcontroller1}) that
\begin{subequations}
\begin{align}
   r (k+1)& = {\underline{A}}{{ r}}(k) +\bar B\bar u(k) +\bar G\vartheta(k), \hfill \\
  {y}(k)& = {{\bar C}}{{ r}}(k). \hfill
\end{align}\label{ceq15}
\end{subequations}
\hspace*{-0.25cm} Based on the policy iteration of (\ref{eq18}) and (\ref{eq19}), (\ref{ceq15}) results in the off-policy Bellman equation for DT systems in the state-feedback form  as
\begin{align}\label{eq27}
  &r^T (k+1)P^{j+1}  r(k+1)-  r^T(k)P^{j+1}  r(k)\nonumber\\
 =& -r^T(k)(\bar Q + (\bar K^j)^T\bar R \bar K^j)r(k)+\vartheta^T(k)\bar G^TP^{j+1}\bar G\vartheta(k) \nonumber\\&
 +(-\bar K^j r(k)+\bar u(k))^T\bar B^TP^{j+1}\bar B(\bar K^j r(k)+\bar u(k))
  \nonumber\\&+2\vartheta^T(k)\bar G^TP^{j+1}\bar Bu(k)+2r^T(k)\underline{A}^TP^{j+1}\bar G\vartheta(k) \nonumber\\&+2r^T(k)\underline{A}^TP^{j+1}\bar B(\bar K^j r(k)+\bar u(k)),
\end{align}
where $\underline{A}^j=\underline{A}-\bar B\bar K^j$.

In (\ref{eq27}), the system {state} $x(k)$ should be known. However, a variety of physical systems, including aircraft and power networks, only allow the measurement of  the system output. Note that, within the output-feedback framework, the unknown term of $r(k)$ (or $x(k)$) prevents us from directly obtaining $\bar K^{k+1}$ from (\ref{eq27}). The next four subsections describe our output-feedback RL for solving \textit{Problem \ref{def2}}.
\subsection{State Reconstruction}
Since the system state $r(k)$ in (\ref{ceq15}) is not available for feedback control, this subsection provides a method to reconstruct $r(k)$ using input-output data from DT systems.

To reconstruct the DT state, we first generate the DT states $\zeta_{\bar u}(k)\in\mathbb R^{n_zr_m}$, $\zeta_y(k)\in\mathbb R^{n_zr_p}$, and $\zeta_\vartheta(k)\in\mathbb R^{n_zr_p}$ through three difference equations with zero initial conditions as
\begin{align}\label{eqvircon7}
{\zeta}_{\bar u}(k+1) & =(I_m\otimes A_\zeta)\zeta_{\bar u}(k)+\bar u(k) \otimes b,\\
{\zeta}_y(k+1) & =(I_p\otimes A_\zeta)\zeta_y(k)+y(k)\otimes b, \label{eqvircon701}\\
{\zeta}_{\vartheta}(k+1) & =(I_{p}\otimes A_\zeta)\zeta_{\vartheta}(k)+\vartheta(k)\otimes b, \label{eqvircon702}
\end{align}
where $A_\zeta$ is a companion matrix  with $-d_{j}$ for $j=1, 2, \cdots, n_z$ at the last row designed to make $A_\zeta$ Schur
and $b=[0, 0, \cdots, 0, 1]^T\in\mathbb R^{n_z}$.

The following result is attained by using the  Luenberger observer for DT systems and extending \cite{Ci2019OPFBRL1}, {Chapters 3.6.1 and 4.5.4 of \cite{tao2003adaptiveBk}}, \cite{Rizvi2018TNNLS}.
\begin{lemma}\label{lemmaLin}
Consider the DT system (\ref{ceq15}). There exists a constant matrix $\bar M\in\mathbb R^{n_z\times r_{\bar\zeta}}$ satisfying
\vspace*{-0.1cm}
\begin{align}\label{eqMM001}
r(k)=\bar M\bar \zeta(k)+\omega(k),
\end{align}
where  $\bar\zeta^T(k)=[\zeta_{\bar u}^T(k), \zeta_y^T(k), \zeta_{\vartheta}^T(k)]^T\in\mathbb R^{r_{\bar\zeta}}$ with $r_{\bar\zeta}={n_zr_m}+{n_zr_p}+{n_zr_p}$   and $\omega(k)=({\underline{A}-\bar L\bar C})^k r(0)$.
\end{lemma}

The matrix $\bar M$
in (\ref{eqMM001}) is termed as a parameterization matrix as in \cite{Ci2019OPFBRL1}. The following result shows that the DT state reconstruction in (\ref{eqMM001}) has a structural property that ${\rm{rank}}(\bar M)$ is linked to the  three controllability matrices $\mathcal{C}(\underline{A}, \bar B)$,
$\mathcal{C}(\underline{A}, \bar L)$, and $\mathcal{C}(\underline{A}, \bar G)$. {This is a necessary and sufficient condition for measuring the rank of the parameterization matrix for DT state reconstruction.}
\begin{lemma}\label{theoremMM}
{\rm{rank}}$(\bar M)=${\rm{rank}}$([\mathcal{C}(\underline{A}, \bar B), \mathcal{C}(\underline{A}, \bar L), \mathcal{C}(\underline{A}, \bar G)])$.
\end{lemma}

\textbf{\textit{Proof:}} See \textit{Appendix \ref{APPtheoremMM}}. \QEDB

This lemma reveals that the  parameterization matrix $\bar{M}$ in the DT version of (\ref{eqMM001}) has the same structural property as that in the CT version of \cite{Ci2019OPFBRL1}. Note that the state reconstruction in (\ref{eqMM001}) is for a tracking problem, which is thus applicable {to} a regulation problem. Based on the property in \textit{Lemma \ref{theoremMM}}, one has the following convergence result for the DT system (\ref{ceq15}).
\begin{theorem}\label{proposition2}
Consider the DT system (\ref{ceq15}) satisfying the controllability condition of $(\underline{A}, \bar B)$ and the observability condition of $(\underline{A}, \bar{C})$,  $r(k)-\bar M\bar \zeta(k)$ asymptotically decays to zero.
\end{theorem}

\textbf{\textit{Proof:}} Under the condition that $(\underline{A}, \bar B)$ is controllable, it follows from \textit{Lemma \ref{theoremMM}} that rank$(\bar M)=n_z$. The matrix $A_\zeta$ in (\ref{eqvircon7})--(\ref{eqvircon702}) is designed to be Schur by choosing appropriate coefficients $d_{j}$ for $j=1, 2, \cdots, n_z$. Thus, the vector $\bar\zeta(k)$ in (\ref{eqMM001}), formed by the difference equations (\ref{eqvircon7})--(\ref{eqvircon702}), is known. In addition, since $(\underline{A}, \bar C)$ is observable, then the eigenvalues of $\underline{A}-\bar L\bar C$ can be designed to be equal to those of $A_\zeta$ through choosing an appropriate matrix $\bar L$. By doing this, $\underline{A}-\bar L\bar C$ is Schur, based on which  $r(k)-\bar M\bar \zeta(k)$ asymptotically decays to zero.  This completes the proof. \QEDB

The reason for seeking rank$(\bar M)=n_z$ is that the uniqueness of the approximate solution to the Bellman equation will depend on it (see \textit{Proof of Lemma \ref{2cclemma2}} in the next subsection).

In the next two subsections, we will focus on how to use the input-output data to learn the optimal control gain matrix after taking the state reconstruction into account.

\subsection{Off-Policy Bellman Equation in Output-Feedback Form}
In this subsection, both $\bar M$ and $\omega(k)$ are introduced to change the state-feedback Bellman equation (\ref{eq27}) into an output-feedback form as

\vspace*{-0.3cm}
\vspace*{-0.1cm}
\begin{align}\label{eq27r1}
  &\bar \zeta^T (k+1)\bar P^{j+1}  \bar \zeta(k+1)-  \bar \zeta^T(k)\bar P^{j+1}\bar \zeta(k)\nonumber\\
 =& -y^T Q y -\bar \zeta^T(k)\bar M^T (\bar K^j)^T\bar R \bar K^j\bar M\bar \zeta(k) \nonumber\\&
 +(-\bar K^j\bar M\bar \zeta(k)+\bar u(k))^T\bar B^TP^{j+1}\bar B(\bar K^j\bar M\bar \zeta(k)+\bar u(k))
  \nonumber\\&+2\vartheta^T(k)\bar G^TP^{j+1}\bar Bu(k)+2\bar \zeta^T(k)\bar M^T\underline{A}^TP^{j+1}\bar G\vartheta(k) \nonumber\\&+2\bar \zeta^T(k)\bar M^T\underline{A}^TP^{j+1}\bar B(\bar K^j \bar M\bar \zeta(k)+\bar u(k))\nonumber\\&+\vartheta^T(k)\bar G^TP^{j+1}\bar G\vartheta(k)+\bar\chi^{j+1}(t),
\end{align}
where
\begin{align}
&\bar\chi^{j+1}(k)\nonumber\\= &-2\omega^T(k+1)P^{j+1}\bar M   \bar\zeta(k+1)-\omega^T(k+1) P^{j+1} \omega(k+1)\nonumber\\
&-2\omega^T(k)P^{j+1}\bar M   \bar\zeta(k)-\omega^T(k) P^{j+1} \omega(k)\nonumber\\
&-2\omega^T(k) (\bar K^j)^T\bar R \bar K^j\bar M   \bar\zeta(k)-\omega^T(k) (\bar K^j)^T\bar R \bar K^j \omega(k)\nonumber\\
&-2\omega^T(k)(\bar K^j)^T\bar B^TP^{j+1}\bar B\bar K^j\bar M   \bar\zeta(k)\nonumber\\
&-\omega^T(k) (\bar K^j)^T\bar B^TP^{j+1}\bar B\bar K^j \omega(k)\nonumber\\
&+2\omega^T(k)\underline{A}^TP^{j+1}\bar G\vartheta(k)
+2\omega^T(k)\underline{A}^TP^{j+1}\bar B\bar K^j\bar M   \bar\zeta(k)\nonumber\\
&+2\omega^T(k) \underline{A}^TP^{j+1}\bar B\bar K^j \omega(k)+2\omega^T(k)\underline{A}^TP^{j+1}\bar Bu(k)\nonumber\\
&+2\bar\zeta^T(k)\bar M^T\underline{A}^TP^{j+1}\bar B\bar K^j \omega(k). \label{eq27r11}
\end{align}

Let
$\mathcal{C}_{\bar\zeta}=[{\rm{vecv}}( \bar\zeta(k_1))-{\rm{vecv}}( \bar\zeta(k_0)), {\rm{vecv}}(\bar\zeta(k_2))-{\rm{vecv}}( \bar\zeta(k_1)), \cdots, {\rm{vecv}}(\bar\zeta(k_f))-{\rm{vecv}}(\bar\zeta(k_{f-1}))]^T$, $\mathcal{D}_{{\bar K}_o^j\bar\zeta}=[{\rm{vecv}}({\bar K}_o^j\bar\zeta(k_0)), {\rm{vecv}}({\bar K}_o^j\bar\zeta(k_1)), \cdots, {\rm{vecv}}({\bar K}_o^j\bar\zeta(k_{f-1}))]^T$, $\mathcal{D}_{\vartheta \bar\zeta}=[\vartheta(k_0)\otimes \bar\zeta(k_0), \vartheta(k_1)\otimes \bar\zeta(k_1),\cdots, \vartheta(k_{f-1})\otimes \bar\zeta(k_{f-1})]^T$, $\mathcal{D}_{\bar\zeta\bar\zeta}=[\bar\zeta(k_0)\otimes \bar\zeta(k_0), \bar\zeta(k_1)\otimes \bar\zeta(k_1), \cdots,$ $ \bar\zeta(k_{f-1})\otimes \bar\zeta(k_{f-1})]^T$, $\mathcal{D}_{\bar u \vartheta}=[\bar u(k_0)\otimes \vartheta(k_0), \bar u(k_1)\otimes \vartheta(k_1),\cdots, \bar u(k_{f-1})\otimes \vartheta(k_{f-1})]^T$,
and $\mathcal{D}_{\bar\chi^{j+1}}=[\bar\chi^{j+1}(t_0),$ $\bar\chi^{j+1}(t_1), \cdots, \bar\chi^{j+1}(t_s)]^T$.  Besides, let $\bar L_P^{j+1}=\bar M^T P^{j+1}\bar M$, $\bar L_1^{j+1}=\bar M^T\underline{A}^TP^{j+1}\bar B$, $\bar L_2^{j+1}=\bar B^TP^{j+1}\bar B$, $L_3^{j+1}=\bar M^T\underline{A}^TP^{j+1}\bar G$, $L_4^{j+1}=\bar G^TP^{j+1}\bar B$, $L_5^{j+1}=\bar G^TP^{j+1}\bar G$, and ${\bar K}_o^j={\bar K}^j\bar M$. Now, (\ref{eq27r1}) is rewritten as
\begin{align}\label{eq28}
{\varrho_o^j} \bar L_{vec}
     ={\nu_o^j}+\mathcal{D}_{\bar\chi^{j+1}},
\end{align}
where
\begin{align}
\bar L_{vec}=&\ [
{\rm{vecs}}^T(\bar L_{P}^{j+1}), 
{\rm{vec}}^T(\bar L_1^{j+1}),
{\rm{vecs}}^T(\bar L_2^{j+1}),\nonumber\\
&\ \ {\rm{vec}}^T(\bar L_3^{j+1}),
{\rm{vec}}^T(\bar L_4^{j+1}),
{\rm{vecs}}^T(\bar L_5^{j+1})
     ]^T,\\
{\varrho_o^j}=&\ [\mathcal{C}_{\bar\zeta}, -2\mathcal{D}_{\bar\zeta \bar\zeta}(I\otimes ({\bar K}_o^j)^T)-2\mathcal{D}_{\bar u \bar\zeta},\nonumber\\& -\mathcal{D}_{\bar u} + \mathcal{D}_{{\bar K}_o^j\bar\zeta}, -2\mathcal{D}_{\vartheta \bar\zeta},-2\mathcal{D}_{\bar u \vartheta}, -\mathcal{D}_\vartheta],\label{eq2801}\\
\nu_o^j=&-\mathcal{D}_{\bar \zeta \bar \zeta}{\rm{vec}}(({\bar K}_o^j)^T\bar R{\bar K}_o^j)-\mathcal{D}_{y y}{\rm{vec}}(Q).\label{eq2802}
\end{align}

From (\ref{eqMM001}), if the initial state of $r(k)$ satisfies $r(0)=0$, then $\omega(k)$ in (\ref{eqMM001}) and  $\mathcal{D}_{\bar\chi^{j+1}}$ in (\ref{eq28}) are zeros. In the next subsection,
we will take the non-zero initial state $r(0)\neq 0$ into account, and seek for approximating (\ref{eq28}).

\subsection{Solution to Output-Feedback Off-Policy Bellman Equation}\label{sectionerror}
In this subsection, we are to reduce the influence from non-zero initials and to give a sufficient condition for approximately solving the off-policy Bellman equation in the output-feedback form.

Here, $\mathcal{D}_{\bar\chi^{j+1}}$ is a nonlinear function of unknown terms $\omega(t)$,  $\bar L_1^{j+1}$, and $\bar L_2^{j+1}$. It thus becomes difficult to obtain an accurate analytical solution from (\ref{eq28}). Instead of directly solving (\ref{eq28}), we turn to computing the following linear equation
\begin{align}\label{eq36r1}
{\varrho_o^j}\hat{\bar L}_{vec}
     ={\nu_o^j},
\end{align}
where $\hat{\bar L}_{vec}=[
{\rm{vecs}}^T(\hat{\bar L}_{P}^{j+1}),
{\rm{vec}}^T(\hat{\bar L}_1^{j+1}),
{\rm{vecs}}^T(\hat{\bar L}_2^{j+1})$,
${\rm{vec}}^T(\hat{\bar L}_3^{j+1}),
{\rm{vec}}^T(\hat{\bar L}_4^{j+1}),
{\rm{vecs}}^T(\hat{\bar L}_5^{j+1})
]^T$, and the notation $\hat{(\cdot)}$ is employed to differ the computed solution in (\ref{eq36r1}) from the analytical one in (\ref{eq28}). We have seen that (\ref{eq36r1}) equals (\ref{eq28}) if $r(0)=0$. In what follows, we focus on handling the non-zero case and give the following result on reducing the solution error between (\ref{eq28}) and (\ref{eq36r1}).

\begin{theorem}\label{theorem3}
Suppose that the non-zero matrix $\bar\varrho^{j}$ in (\ref{eq2801}) is collected over the time interval $[k_0, k_f]$.  If the starting time for the data collection $k_0$ is sufficiently large, then the computed solution from (\ref{eq36r1}) can be
 considered as an approximate solution of (\ref{eq28}) with the solution error being sufficiently small.
\end{theorem}

\textbf{\textit{Proof:}}  Consider two difference equations
\begin{align}
 v(s+1)=&v(s)-\varepsilon{(\varrho_o^j)}^T({\varrho_o^j}v(s)-\nu_o^j-\mathcal{D}_{\bar\chi^{j+1}}),\label{eq36r2}\\
\hat v(s+1)=&\hat v(s)-\varepsilon{(\varrho_o^j)}^T({\varrho_o^j}\hat v(s)-\nu_o^j),\label{eq36r3}
\end{align}
to solve (\ref{eq28}) and (\ref{eq36r1})
with $v(0)=\hat v(0)=0$ and the constant $\varepsilon$ satisfying
\begin{align}\label{eq11111}
0<\varepsilon<{2\rho^{-1}({(\varrho_o^j)}^T({\varrho_o^j}))}.
\end{align}
 Bringing the algorithmic time $s$ into (\ref{eq36r2}) and (\ref{eq36r3}) is to distinguish it from the system evolution time $k$ used in (\ref{eqsystem7})--(\ref{eqsystem72}).

 It follows from the singular value decomposition that there exist matrices $W$ and $\varrho_w^j$ with $W^TW=WW^T=I$ and $(\varrho_w^j)^T\varrho_w^j>O$ such that\cite{Ci2019OPFBRL1,liu2018DTsolver}
 $\varrho_o^jW =[\varrho_w^j\ O]$.
Define \begin{align}\label{cceq36r503}
\bar v(s) = W^Tv(s)=[{\bar v}_{1}^T(s), {\bar v}_{2}^T(s)]^T.
 \end{align}
 From
 (\ref{cceq36r503}), (\ref{eq36r2}) is rewritten into
\begin{align}\label{eq36r503}
{\bar v}_{1}(s+1)=&{\bar v}_{1}(s)-\varepsilon{(\varrho_w^j)}^T\varrho_w^j{\bar v}_{1}(s)\nonumber\\&+\varepsilon{(\varrho_w^j)}^T\nu_o^j+\varepsilon{(\varrho_w^j)}^T\mathcal{D}_{\bar\chi^{j+1}},\\
{\bar v}_{2}(s+1)=&{\bar v}_{2}(s),\label{eq36r504}
\end{align}
where ${\bar v}_2(s)$ is zero for any $s$. Hence, (\ref{eq36r503}) is changed into
 \begin{align}\label{cceq36r504}
&||{\bar v}_{1}(s+1)||\nonumber\\
\leq& \sqrt{|\lambda_{\max}(I-\varepsilon{(\varrho_w^j)}^T\varrho_w^j)|}||{\bar v}_{1}(s)||+||\varepsilon{(\varrho_w^j)}^T\nu_o^j||\nonumber\\&+||\varepsilon{(\varrho_w^j)}^T\mathcal{D}_{\bar\chi^{j+1}}||.
 \end{align}
 From (\ref{eq11111}), one has that  $0<\rho(\varepsilon{(\varrho_w^j)}^T\varrho_w^j)<2$ such that $\rho(I-\varepsilon{(\varrho_w^j)}^T\varrho_w^j)<1$. Therefore, given the constant $\varepsilon$ and the matrix ${(\varrho_w^j)}^T\varrho_w^j$, there exists a constant $\epsilon$ so that $|\lambda_{\max}(I-\varepsilon{(\varrho_w^j)}^T\varrho_w^j)|<\epsilon<1$. As a result, one changes (\ref{cceq36r504}) into
 \begin{align}\label{cceq36r505}
||{\bar v}_{1}(s+1)|| \leq& \epsilon_1 ||{\bar v}_{1}(s)||+||\varepsilon{(\varrho_w^j)}^T\nu_o^j||\nonumber\\&+||\varepsilon{(\varrho_w^j)}^T\mathcal{D}_{\bar\chi^{j+1}}||,
\end{align}
where $\epsilon_1=\sqrt{|\lambda_{\max}(I-\varepsilon{(\varrho_w^j)}^T\varrho_w^j)|}<\sqrt{\epsilon}<1$. It is noted that the term $||\varepsilon{(\varrho_w^j)}^T\nu_o^j||$ in (\ref{cceq36r505}) is bounded.

Now, we focus on the boundedness of $\mathcal{D}_{\bar\chi^{j+1}}$ in (\ref{cceq36r505}).
Let $\alpha_{i}$ for $i=0,1,2,3,4$ and $\epsilon_l$ for $l=1,2$ be certain positive constants. Since $\underline{A}-\bar L\bar C$ is Schur, $0<|\lambda_{\max}(\underline{A}-\bar L\bar C)|<1$ holds.
From (\ref{eq27r11}), one has
\begin{align}\label{eq36r6}
||\mathcal{D}_{\bar\chi^{j+1}}|| \leq&\  ||{\bar v(s)}||\alpha_{2}|\lambda_{\max}(\underline{A}-\bar L\bar C)|^{k_0}.
\end{align}
 Substituting (\ref{eq36r6}) into (\ref{cceq36r505}) yields
\begin{align}
||{\bar v}_{1}(s+1)||
 \leq\ &(\epsilon_1+\alpha_{3}|\lambda_{\max}(\underline{A}-\bar L\bar C)|^{k_0})||{\bar v}_{1}(s)||+\bar b^j,\nonumber
\end{align}
where $\bar b^j$ is defined as an upper bound of $||\varepsilon{(\varrho_w^j)}^T\mathcal{D}_{\bar\chi^{j+1}}||$. Thus, with sufficiently large $k_0$, $\alpha_{3}|\lambda_{\max}(\underline{A}-\bar L\bar C)|^{k_0}$  is thus sufficiently small. There exists a positive constant $\epsilon_2$ satisfying $0<\epsilon_1<\epsilon_2<1$ and $\epsilon_2-\epsilon_1>\alpha_{3}|\lambda_{\max}(\underline{A}-\bar L\bar C)|^{k_0}$ such that
\begin{align}\label{ccceq36r4}
||{\bar v}_{1}(s+1)||
          \leq& \epsilon_2 ||{\bar v}_{1}(s)||+\bar b^j,
\end{align}
which leads to
\begin{align}\label{cc1ceq36r4}
||{\bar v}_{1}(s)||\leq &\Big(||{\bar v}_{1}(0)||-\frac{\bar b^j}{1-\epsilon_2}\Big)\epsilon_2^s+\frac{\bar b^j}{1-\epsilon_2}.
\end{align}
Thus, $\bar v(s)$, $v(s)$, and $\mathcal{D}_{\bar\chi^{j+1}}$ are bounded for any $s$.

With the boundedness of $v(s)$, we now prove the convergence of the solution error between (\ref{eq28}) and  (\ref{eq36r1}) in the remaining analysis.
Similar to (\ref{cceq36r503}), we define
\begin{align}\label{eq36r50001}
\hat{\bar v}(s)=W^T \hat v(s)=[\hat{\bar v}_1^T(s), \hat{\bar v}_2(s)]^T,
\end{align}
based on which one has
\begin{align}\label{eq36r505}
\hat{\bar v}_{1}(s+1) &=\hat{\bar v}_{1}(s)-{(\varrho_w^j)}^T{\varrho_w^j}\hat{\bar v}_{1}(s)+{(\varrho_w^j)}^T\nu_o^j,\\
\hat{\bar v}_{2}(s+1) &=\hat{\bar v}_{2}(s),\label{eq36r506}
\end{align}
where $\hat{\bar v}(s)$ is zero. Thus, from (\ref{eq36r506}), $\hat{\bar v}_{2}(s)=0$ holds. Let the error be $\tilde{\bar v}_{1}(s)=\bar v_1(s)-\hat{\bar v}_{1}(s)$. From (\ref{eq36r503}) and (\ref{eq36r505}), one has
\begin{align}\label{eq36r4}
||\tilde{\bar v}_{1}(s+1)||
 \leq\ &\epsilon_1||\tilde{\bar v}_{1}(s)||+\alpha_{4}|\lambda_{\max}(\underline{A}-\bar L\bar C)|^{k_0},
\end{align}
where the boundedness of $\mathcal{D}_{\bar\chi^{j+1}}$ is employed.
By (\ref{eq36r4}), one has
\begin{align}\label{11eq36r4}
||\tilde{\bar v}_{1}(s)||\leq\ &\Big(||{\bar v}_{e_1}(0)||-\frac{\alpha_{4}|\lambda_{\max}(\underline{A}-\bar L\bar C)|^{k_0}}{1-\epsilon_1}\Big)|\epsilon_1|^s
\nonumber\\&+\frac{\alpha_{4}|\lambda_{\max}(\underline{A}-\bar L\bar C)|^{k_0}}{1-\epsilon_1},
\end{align}
where the constant $\epsilon_1$ has been defined in (\ref{cceq36r505}) satisfying $0<\epsilon_1<1$ and the term $\alpha_{4}|\lambda_{\max}(\underline{A}-\bar L\bar C)|^{k_0}$ has been tuned sufficiently small under sufficiently large $k_0$. From (\ref{cceq36r503}) and (\ref{eq36r50001}),
one changes (\ref{11eq36r4}) into
\begin{align}\label{eq36r50002}
\lim_{s\rightarrow\infty}||\hat v(s)-v(s)|| \leq & \lim_{s\rightarrow\infty}||W||||\hat{\bar v}(s)-\bar v(s)||\nonumber
\\=& \lim_{s\rightarrow\infty}||\hat{\bar v}(s)-\bar v(s)||\nonumber\\
=&\ \frac{\alpha_{4}|\lambda_{\max}(\underline{A}-\bar L\bar C)|^{k_0}}{1-\epsilon_1},
\end{align}
which implies that the solution of (\ref{eq28}) converges to that of (\ref{eq36r1}) with the error being sufficiently small by increasing $k_0$. Therefore, the proof is completed.
\QEDB
\begin{remark}\label{lemma8}
\textit{Theorem \ref{theorem3}} shows that the solution error between (\ref{eq28}) and (\ref{eq36r1}) can be made smaller by choosing a larger starting time $k_0$ for the data collection. This prompts us to introduce an additional \textit{Model-Free Pre-Collection Phase}. The necessity of the additional phase roots in the requirement of the convergence in the DT state reconstruction and the idea is inspired by the CT work of \cite{Ci2019OPFBRL1}. Details on how to coordinately execute the DT off-policy learning will be presented in \textit{Algorithm 1}. \QEDD
\end{remark}
\begin{remark}
In contrast to \cite{Ci2019OPFBRL1} wherein differential equations are used to solve the linear equation, \textit{Theorem \ref{theorem3}} uses difference equations. In addition, we find that our result in \textit{Theorem \ref{theorem3}} and that in \cite{huangDT2016TAC} are dual to each other in {a} certain sense. To be specific, \textit{Theorem \ref{theorem3}} shows that a class of unknown nonlinear equations are approximately solved by a known linear difference equation, while \cite{huangDT2016TAC} shows that an unknown linear difference equation is approximately solved by a class of nonlinear equations. \QEDD
\end{remark}

The following result shows that how much data should we collect to achieve the optimal output tracking control within the output-feedback RL framework.
\begin{lemma}\label{2cclemma2}
The off-policy Bellman equation in the output-feedback form (\ref{eq36r1}) over the time interval $[k_0, k_f]$,  $\hat{\bar L}_i^{j+1}$, for $i=1,2,\ldots, 5$, can be uniquely solved, if
\begin{enumerate}
  \item  the collected input-output data  at the time $k_f$ satisfy
  \begin{align}\label{eq29}
& {\rm{rank}}([\mathcal{D}_{\bar\zeta \bar\zeta}, \mathcal{D}_{\bar u \bar\zeta}, \mathcal{D}_{\bar u}, \mathcal{D}_{ \vartheta \bar\zeta}, \mathcal{D}_{\bar u \vartheta}, \mathcal{D}_\vartheta])\nonumber\\&=\frac{1}{2}(n_{\bar\zeta}(n_{\bar\zeta}+1)+r_m(r_m+1)+r_p(r_p+1))\nonumber\\&\quad\ +n_{\bar\zeta}r_m+r_pn_{\bar\zeta}+r_pr_m;
\end{align}
  \item rank$(\bar M)=n_z$;
  \item $r(0)=0$.
\end{enumerate}
\end{lemma}

\textbf{\textit{Proof:}}
See \textit{Appendix \ref{lemmaC}}. \QEDB

In \textit{Lemma \ref{2cclemma2}}, the condition $r(0)=0$ is required. We now extend it to the condition $r(0)\neq 0$ by choosing the starting time $k_0$ to be sufficiently large. This is because the solution $\hat{\bar L}_i^{j+1}$, for $i=1,2,\ldots, 5$ in (\ref{eq36r1}) under the condition $r(0)\neq 0$ converges to the unique solution in \textit{Lemma \ref{2cclemma2}} under the condition $r(0)=0$ after recalling the result in \textit{Theorem \ref{theorem3}}.

\subsection{Optimal Output Tracking Design via Output-Feedback RL}
Our output-feedback off-policy learning algorithm is summarized in \textit{Algorithm 1}, where the successive error of the control matrix $\| {\hat{\bar K}_o^{j+1} - \hat{\bar K}_o^{j}} \|$ is used for the stopping indicator since $\hat{\bar K}_o^{j+1}$ is solved uniquely under conditions in \textit{Lemma \ref{2cclemma2}}.

Based on the learned optimal control gain matrix  $\hat{\bar K}_o^*$ from \textit{Algorithm 1}, we achieve optimal output tracking control of DT systems via output feedback as below.

\begin{algorithm}[t]
  \caption{\textit{Output-Feedback Off-Policy RL for Optimal Output Tracking Control of DT Systems}}
  \begin{algorithmic}[1]
    \State\textbf{Initialize:} {Let $j=0$ and ${\bar K}_o^j$ be a stabilizing gain. Let the pair $(F, G)$ be an $r_p$-copy internal model of $S$.}
\State \textbf{Model-Free Pre-Collection Phase:}  From (\ref{eqvircon7})--(\ref{eqvircon702}), compute ${\zeta}_u$, ${\zeta}_y$, and $\zeta_\vartheta$ over the time interval $[0, k_0)$, where $k_0$ is set sufficiently large.
    \State \textbf{Data-Collection Phase:}
   Apply the behavior policy $\bar u(k)$ in (\ref{ceqcontroller1}) with $(F, T)$ being observable to the system (\ref{eqsystem7}) with $r_p\geq r_m$ over the time interval $[k_0, k_f]$  for collecting the input-output data $\mathcal{D}_{\bar\zeta \bar\zeta}$, $\mathcal{D}_{\bar u \bar\zeta}$, $\mathcal{D}_{\bar u}$, $\mathcal{D}_{ \vartheta \bar\zeta}$, $\mathcal{D}_{\bar u \vartheta}$, and $\mathcal{D}_\vartheta$.
     \If{(\ref{eq29}) holds}
     \While{the stopping indicator $\| {\hat{\bar K}_o^{j+1} - \hat{\bar K}_o^{j}} \| \leqslant {\varepsilon }$ is not satisfied with  ${\varepsilon } > 0$ being a small constant}
     \State Solve (\ref{eq36r1}) and update the feedback gain as
                              \begin{align}\label{22eqalg2}
                              \hat{\bar K}_o^{j+1}=(\bar R + \hat{\bar L}_2^{j+1})^{-1}(\hat{\bar L}_1^{j+1})^T.
                              \end{align}
     \EndWhile
          \EndIf
                                \State \textbf{Optimal Control Phase:}
                       The learned optimal control gain matrix $\hat{\bar K}_o^*$
                       is given as $\hat{\bar K}_o^*=\hat{\bar K}_o^{j+1}$.
  \end{algorithmic}
\end{algorithm}

\begin{theorem}\label{gtheorem3}
Let the DT system (\ref{eqsystem7})--(\ref{eqsystem72}) satisfy \textit{Assumptions \ref{assumption3}--\ref{assumption4}} and the output-feedback adaptive optimal output tracking DT controller be designed as
\begin{align}
u(k)=&-\hat{\bar K}_o^*\bar \zeta(k)-T{\bar z(k)},\label{eqcontroller}
\end{align}
where $\bar \zeta$ is given in (\ref{eqMM001}), $\bar z(k)$ is given in (\ref{ceqcontroller1}) with $\vartheta(k)$ being assigned as $y_d(k)$, and  $\hat{\bar K}_o^*$ is the {learnt} optimal control gain matrix by \textit{Algorithm 1}. Therefore, the state-oriented tracking
optimization problem, \textit{Problem \ref{def2}}, is solved with the output tracking error $y_e(k)$ in (\ref{eqsystem8}) decaying to zero, asymptotically.
\end{theorem}

\textbf{\textit{Proof:}}
{We first show that the learnt matrix  $\hat{\bar K}_o^{j+1}$ in  (\ref{22eqalg2}) converges to the ideal matrix $\bar K^*\bar M$ with $\bar K^*$ and $\bar M$ being defined in (\ref{eq1601}) and (\ref{eqMM001}), respectively. With the condition in \textit{Theorem \ref{theorem3}} satisfied, $\hat{\bar L}_1^{j+1}$ and $\hat{\bar L}_2^{j+1}$ in (\ref{eq36r1}) can be made to converge to ${\bar L}_1^{j+1}$ and ${\bar L}_2^{j+1}$ in (\ref{eq28}), respectively. This, together with the uniqueness in \textit{Lemma \ref{2cclemma2}}, leads that the learned matrix $(\bar R + \hat{\bar L}_2^{j+1})^{-1}(\hat{\bar L}_1^{j+1})^T$ in (\ref{22eqalg2}) uniquely converges to $\bar K^{j+1}\bar M$, where $\bar K^{j+1}$ is defined in (\ref{eq19}). It follows from \textit{Lemma \ref{lemma1}} that $\bar K^{j+1}$ converges to $\bar K^*$ as the integer $j$ gets larger. Note that $\bar K^*$ is the unique solution satisfying the \textit{Problem \ref{def2}} under the controllability of $(\underline{A}, \bar B)$ and the observability of $(\underline{A}, \bar{C})$. Therefore, if $\hat{\bar K}_o^{j+1}$ converges, then the unique solution  $\hat{\bar K}^{j+1}$ from (\ref{eq36r1}) converges to $\bar K^*\bar M$.  It reveals that the ideal matrix $\bar K^*\bar M$ is learned by the matrix $\hat{\bar K}_o^*$ from \textit{Algorithm 1}.} We then show the convergence of the output tracking error $y_e(k)$.
 With the learned optimal control gain matrix $\hat{\bar K}_o^*$,  the closed-loop system becomes
$e(k+1) = {{\bar A}_o^*}e(k)+ \bar B{(\underline{A}-\bar L\bar C)^k}r(0)$, where $e(k)$ is given in (\ref{eqerror1pre}) and both
${\bar A_o^*}{\text{ = }}\underline{A}-\bar B\hat{\bar K}_o^*$ and $\underline{A}-\bar L\bar C$ are Schur. Considering that $\lim_{k\rightarrow\infty}\bar B{(\underline{A}-\bar L\bar C)^k}r(0)=0$, one obtains that $\lim_{k\rightarrow\infty}e(k)=0$\cite{huangDT2016TAC}, based on which the output tracking error satisfies $\lim_{k\rightarrow\infty}y_e(k)=0$ from (\ref{eq1901pre02}).  This completes the proof. \QEDB
 \begin{remark}
The RL-based controller in (\ref{eqcontroller})  is robust to system uncertainties, which corresponds to the linear robust output regulation in the literature such as \cite{huang2004nonlinear}. That is, after the learning is completed, the proposed controllers are robust to some model uncertainties in the system dynamics matrices $A+\Delta A$, $B+\Delta B$, and $C+\Delta C$, where $A$, $B$, and $C$ denote the nominal part of the plant; $\Delta A$, $\Delta B$, and $\Delta C$ represent the model uncertainties.
\QEDD
 \end{remark}

\section{Conclusion}\label{Sect6}
This paper investigated the output-feedback optimal output tracking problem for DT systems with unknown dynamics using the off-policy RL and {robust output regulation theory}. We have formulated the output tracking optimization problem based on the newly proposed dynamical DT controller in contrast to the standard DT controller by linear output regulation theory. We have shown that, by making use of the collected data {along with} the controlled system and the reference output, we are able to approximate the optimal output-feedback controller. We have studied the parameterization matrix and re-expression
error so that the learned optimal controller has a satisfactory performance.

\appendices
\section{Proof of \textit{Lemma \ref{theoremMM}}}\label{APPtheoremMM}
Let us consider a standard Luenberger observer as
 \begin{align}\label{eqobs}
{\hat r}(k+1)=&\underline{A}\hat r(k) + \bar B \bar u(k)  + \bar G\vartheta(k)\nonumber\\
&+\bar L (y(k)-\bar C\hat r(k))
\end{align}
with $\bar L$ being a $n_z\times r_p$ matrix
so that $\hat r(k)-r(k)$ decays to zero with $r(k)$ {given} by (\ref{ceq15}).
Similar to \cite{Rizvi2018TNNLS},  (\ref{eqobs}) is rewritten as 
\begin{align}\label{eqvircon4}
\hat r(k)=&\ G_u(z)[\bar u(k)]+G_y(z)[y(k)]\nonumber\\&+G_\vartheta(z)[\vartheta(k)]+{(\mathcal{A}_o)^k}\hat r(0),
\end{align}
where $\mathcal{A}_o=\underline{A}-\bar L\bar C$, $G_{\bar u}(z)[\bar u(k)]$ is a time-domain DT {signal} represented by the frequency-domain representation  $G_{\bar u}(z)=\left[ {\begin{array}{*{20}{cc}}
  G_{1,1}^{\bar u}(z) &\cdots  &G_{1,r_m}^{\bar u}(z) \\[-0.2cm]
  \vdots &\ddots  &\vdots \\[-0.2cm]
    G_{n_z,1}^{\bar u}(z) &\cdots  &G_{n_z,r_m}^{\bar u}(z) \\
\end{array}} \right]\in\mathbb R^{n_z\times r_m}$. The entry $G_{i,j}^{\bar u}(z)$ at the $i$th row and $j$th column of $G_{\bar u}(z)$ is extended to
\begin{align}\label{eqvircon5}
G_{i,j}^{\bar u}(z)=&\frac{1}{\det(zI-\mathcal{A}_o)}\Big[g_{i,j,n_z-1}^{\bar u}z^{n_z-1}+g_{i,j,n_z-2}^{\bar u}z^{n_z-2}\nonumber\\&\quad\quad\quad+\cdots+g_{i,j,1}^{\bar u}z+g_{i,j,0}^{\bar u}\Big].
\end{align}
Here, $G_{\bar u}(z)[\bar u(k)]$ is interpreted as
\begin{align}\label{eqvircon6}
&G_{\bar u}(z)[\bar u(k)]=\nonumber\\&\hspace*{-0.3cm}\underbrace{\left[ {\begin{array}{*{20}{cc}}
  g_{1,1,0}^{\bar u} &\cdots & g_{1,r_m,n_z-1}^{\bar u}& g_{1,r_m+1,0}^{\bar u} &\cdots  &g_{1,n_zr_m,n_z-1}^{\bar u} \\[-0.2cm]
  \vdots &\ddots  &\vdots &\vdots &\ddots  &\vdots \\[-0.2cm]
   g_{n_z,1,0}^{\bar u} &\cdots & g_{n_z,r_m,n_z-1}^{\bar u}& g_{1,r_m+1,0}^{\bar u} &\cdots  &g_{n_z,n_zr_m,n_z-1}^{\bar u} \\
\end{array}} \right]}_{M_{\bar u}}\nonumber\\
&\quad\quad\times\underbrace{[u(k)]\otimes \left[1,  z,  \cdots,  z^{n_z-1} \right]^T\frac{1}{\det(zI-\mathcal{A}_o)}}_{\zeta_{\bar u}(k)}\nonumber\\
&\triangleq M_{\bar u}\zeta_{\bar u}(k).
\end{align}
Similar definitions of $M_y\zeta_y(k)$ and $M_\vartheta\zeta_\vartheta(k)$ apply to $G_y(z)[y(k)]$ and $G_\vartheta(z)[\vartheta(k)]$, respectively. Thus, one obtains that
\begin{align}\label{eqMM}
\hat r(k)&=\bar M\bar \zeta(k)+{(\mathcal{A}_o)^k}\hat r(0).
\end{align}

From (\ref{eqobs}) and the z-transform operator, one has
\begin{align}
\hat r(k)=&{(\underline{A}-\bar L\bar C)^k}\hat r(0)+ (zI-\mathcal{A}_o)^{-1}\bar B[\bar u(k)]\nonumber\\&+(zI-\mathcal{A}_o)^{-1}\bar L[y(k)]+(zI-\mathcal{A}_o)^{-1}\bar G[\vartheta(k)],
\nonumber\end{align}
where
\begin{align}
(zI-\mathcal{A}_o)^{-1}=&\frac{{\rm{adj}}(zI-\mathcal{A}_o)}{\det(zI-\mathcal{A}_o)},\label{eqvircon041}
\end{align}
\begin{align}
\det(zI-\mathcal{A}_o)=&z^{n_z}+d_1z^{{n_z}-1}+d_2z^{{n_z}-2}+\cdots\nonumber\\&+d_{{n_z}-1}z+d_{n_z},\label{eqvircon042}\\
{\rm{adj}}(zI-\mathcal{A}_o)=&B_0z^{n_z-1}+B_1z^{{n_z}-2}+\cdots\nonumber\\&+B_{{n_z}-2}z+B_{{n_z}-1}. \label{eqvircon043}
\end{align}
From (\ref{eqvircon041}) and (\ref{eqvircon043}), one has
\begin{align}\label{eqvircon044}
&\det(zI-\mathcal{A}_o)I\nonumber\\=&\ B_0z^{n_z}+(B_1-B_0\mathcal{A}_o)z^{{n_z}-1}+(B_2-B_1\mathcal{A}_o)z^{{n_z}-2}\nonumber\\&
               +\cdots+(B_{{n_z}-1}-B_{{n_z}-2}\mathcal{A}_o)z-B_{{n_z}-1}A.
\end{align}
Using (\ref{eqvircon042}) and (\ref{eqvircon044}),  the following equations $B_0=I$ and $B_{i+1}=B_i\mathcal{A}_o+d_{i+1}I$ hold for  $i=0, 1, 2, \cdots, n_z-2$.
 For $i=0$, one has $B_{1}=B_0\mathcal{A}_o+d_{1}I=\mathcal{A}_o+d_{1}I=\underline{A}-\bar L\bar C+d_{1}I$. This leads to
 \begin{align}\label{eqvircon047}
&{\rm{rank}}([[B_0, B_1]\bar B, [B_0, B_1]\bar L, [B_0, B_1]\bar G])
\nonumber\\=\ &{\rm{rank}}([[B_0, \underline{A}-\bar L\bar C+d_{1}I]\bar B, [B_0, \underline{A}-\bar L\bar C+d_{1}I]\bar L,\nonumber\\&\quad\quad [B_0, \underline{A}-\bar L\bar C+d_{1}I]\bar G])
\nonumber\\=\ &{\rm{rank}}([[\bar B, (\underline{A}-\bar L\bar C)\bar B+d_{1}\bar B], [\bar L, (\underline{A}-\bar L\bar C)\bar L+d_{1}\bar L],\nonumber\\&\quad\quad [\bar G, (\underline{A}-\bar L\bar C)\bar G+d_{1}\bar G]])
\nonumber\\=\ &{\rm{rank}}([\bar B, \underline{A}\bar B, \bar L, \underline{A}\bar L, \bar G, \underline{A}\bar G]),
\end{align}
where the last equation is obtained using the fact that column operations do not change the rank of a matrix.
Based on analysis in (\ref{eqvircon047}) with $i=0$, proceeding the order $i$ to be higher one by one, one has
\begin{align}\label{eqvircon052}
&{\rm{rank}}([\mathcal{B}(\bar B), \mathcal{B}(\bar L), \mathcal{B}(\bar G)])\nonumber\\=\ &{{\rm{rank}}}([\mathcal{C}(\underline{A}, \bar B), \mathcal{C}(\underline{A}, \bar L), \mathcal{C}(\underline{A}, \bar G)])
\end{align}
with $\mathcal{B}(\bar B)=[B_0, B_1, \cdots, B_{n_z-1}]\bar B$, $\mathcal{B}(\bar L)=[B_0, B_1, \cdots$, $B_{n_z-1}]\bar L$, and $\mathcal{B}(\bar G)=[B_0, B_1, \cdots, B_{n_z-1}]\bar G$.

Based on (\ref{eqvircon052}), we next clarify that ${\rm{rank}}(\bar M)={\rm{rank}}([\mathcal{B}(\bar B), \mathcal{B}(\bar L), \mathcal{B}(\bar G)])$ with $\bar M$ being given in (\ref{eqMM}).
Based on the defintions in (\ref{eqvircon6}) and (\ref{eqvircon043}), ${\rm{rank}}(M_{\bar u})={\rm{rank}}(\mathcal{B}(\bar B))$. It leads to  ${\rm{rank}}(M_{y})={\rm{rank}}(\mathcal{B}(\bar L))$ and ${\rm{rank}}(M_{\vartheta})={\rm{rank}}(\mathcal{B}(\bar G))$. Hence, one has  ${\rm{rank}}(\bar M)={\rm{rank}}([\mathcal{B}(\bar B), \mathcal{B}(\bar L), \mathcal{B}(\bar G)])$. This, together with (\ref{eqvircon052}) leads to rank$(\bar M)=$rank$([\mathcal{C}(\underline{A}, \bar B), \mathcal{C}(\underline{A}, \bar L), \mathcal{C}(\underline{A}, \bar G)])$. Therefore,  the proof is completed. \QEDB

\section{Proof of \textit{Lemma \ref{2cclemma2}}}\label{lemmaC}

In order to show the uniqueness of $\hat{\bar L}_i^{j+1}$, for $i=1,2,\ldots, 5$, in (\ref{eq36r1}), it is equivalent to proving that  
 \begin{align}\label{22eq30}
0=&{\varrho_o^j}\bar\Xi^v,
\end{align}
with $\bar\Xi^v=[\bar W^v,\bar Y_1^v,\bar Y_2^v,\bar Y_3^v,\bar Y_4^v,\bar Y_5^v]^T$  has a unique zero solution $[\bar Y_1^v,\bar Y_2^v,\bar Y_3^v,\bar Y_4^v,\bar Y_5^v]^T=0$, where $\bar W^v={\rm{vecs}}(\bar W^m)$, $\bar  Y_1^v={\rm{vec}}(\bar Y_1^m)$, $\bar Y_2^v={\rm{vecs}}(\bar Y_2^m)$, $\bar Y_3^v={\rm{vec}}(\bar Y_3^m)$, $\bar Y_4^v={\rm{vec}}(\bar Y_4^m)$, and $\bar Y_5^v={\rm{vecs}}(\bar Y_5^m)$ with $\bar W^m=(\bar W^m)^T$, $\bar Y_2^m=(\bar Y_2^m)^T$, and $\bar Y_5^m=(\bar Y_5^m)^T$.  Define
\begin{align}\label{22eqcc1}
\bar Z^m=(\bar M \bar M^T)^{-1}\bar M \bar W^m \bar M^T(\bar M\bar M^T)^{-1},
\end{align}
where the property of rank$(\bar M)=n_z$ is employed.  {Under the condition (3) of \textit{Lemma \ref{2cclemma2}}, the approximation errors  $\bar\chi^{j+1}(t)$ in (\ref{eq27r11}) are zeros. Thus, it follows from (\ref{eq27r1}) and (\ref{22eqcc1}) that (\ref{22eq30}) leads to
\begin{align}\label{22cceq31}
0&=\ \mathcal{D}_{\bar\zeta}{\rm{vecs}}(\bar\kappa_P)+2\mathcal{D}_{\bar u \bar\zeta}{\rm{vec}}(\bar\kappa_1)+\mathcal{D}_{\bar u} {\rm{vecs}}(\bar\kappa_2)\nonumber\\&+2\mathcal{D}_{ \vartheta \bar\zeta}{\rm{vec}}(\bar\kappa_3)+2\mathcal{D}_{\bar u \vartheta}{\rm{vec}}(\bar\kappa_4)+ \mathcal{D}_\vartheta {\rm{vecs}}(\bar\kappa_5),
\end{align}}
where
$\bar\kappa_{P}=\bar M^T[(\underline{A}^j)^T\bar Z^m\underline{A}^j-\bar Z^m]\bar M+ ({\bar K}_o^j)^T(\bar B^T\bar Z^m\bar B-\bar Y_2^m){\bar K}_o^j+(\underline{A}^T\bar Z^m\bar B-\bar Y_1^m){\bar K}_o^j+({\bar K}_o^j)^T(\underline{A}^T\bar Z^m\bar B-Y_1^m)^T$,
$\bar\kappa_{1}=\underline{A}^T\bar Z^m\bar B-\bar Y_1^m$,
$\bar\kappa_2=\bar B^T\bar Z^m\bar B-\bar Y_2^m$
$\bar\kappa_{3}=\underline{A}^T\bar Z^m\bar G-\bar Y_3^m$,
$\bar\kappa_4=\bar G^T\bar Z^m\bar B-\bar Y_4^m$, and
$\bar\kappa_5=\bar G^T\bar Z^m\bar G-\bar Y_5^m$.

The matrix $\big[\mathcal{D}_{\bar\zeta}, 2\mathcal{D}_{\bar u \bar\zeta}, \mathcal{D}_{\bar u}, 2\mathcal{D}_{\vartheta \bar\zeta}, 2\mathcal{D}_{\bar u \vartheta}, \mathcal{D}_\vartheta\big]$ is full column rank if (\ref{eq29}) holds. Thus, the solution to (\ref{22cceq31}) is uniquely obtained as
\begin{align}\label{22eqxxxcc1}
(&{\rm{vecs}}^T(\bar\kappa_{P}),
{\rm{vec}}^T(\bar\kappa_{1}),
{\rm{vecs}}^T(\bar\kappa_{2}),
\nonumber\\&
{\rm{vec}}^T(\bar\kappa_{3}),
{\rm{vec}}^T(\bar\kappa_{4}),
{\rm{vecs}}^T(\bar\kappa_{5}))^T=0.
\end{align}
 Recalling rank$(\bar M)=n_z$, one further rewrites $\bar\kappa_{P}$ in (\ref{22cceq31}) as
$(\underline{A}^j)^T\bar Z^m\underline{A}^j-\bar Z^m=0$,
 where $\underline{A}^j$ is Schur. Therefore,  $\bar Z^m$ must be zeros, based on which $\bar Y_i^v$ in (\ref{22eq30}), for $i=1,2,\ldots, 5$, are also zeros. {This implies that $\hat{\bar L}_i^{j+1}$ in (\ref{eq36r1}), for $i=1,2,\ldots, 5$,  are unique. {Note that since the non-square matrix $\bar M$ in (\ref{22eqcc1}) is only full row rank, the zero solution of $\bar Z^m$ does not ensure the zero solution of $\bar W^m$ so that $\hat{\bar L}_P^{j+1}$ may not be unqiue.} This completes the proof. 
\QEDB

\section*{Acknowledgement}
The authors are thankful to Prof. Zongli Lin for his helpful comments and suggestions. 
\ifCLASSOPTIONcaptionsoff
  \newpage
\fi



%
\bibliographystyle{IEEEtran}

\end{document}